\documentclass[a4paper,reqno]{amsart}

\usepackage[utf8]{inputenc}
\usepackage[english]{babel}
\usepackage{amsmath,amssymb,amsthm,exscale,amsopn,mathtools,color,tikz-cd,cite,fullpage}
\usepackage[colorlinks=true,pdftex,unicode=true,linktocpage,bookmarksopen,hypertexnames=false]{hyperref}

\newtheorem{lem}{Lemma}[section]
\newtheorem{prop}[lem]{Proposition}
\newtheorem{thm}[lem]{Theorem}
\newtheorem{cor}[lem]{Corollary}

\theoremstyle{definition}
\newtheorem{ex}[lem]{Example}
\newtheorem{rem}[lem]{Remark}

\newtheorem{notation}[lem]{Notation}

\DeclareMathOperator{\Aut}{Aut}
\DeclareMathOperator{\ch}{char}
\DeclareMathOperator{\clK}{clKdim}
\DeclareMathOperator{\End}{End}
\DeclareMathOperator{\GK}{GKdim}
\DeclareMathOperator{\gr}{gr}
\DeclareMathOperator{\id}{id}
\DeclareMathOperator{\Img}{Im}

\DeclareMathOperator{\Map}{Map}
\DeclareMathOperator{\rk}{rk}
\DeclareMathOperator{\Sym}{Sym}

\newcommand{\free}[1]{\langle #1 \rangle}

\newcommand{\mb}{\mathbb}
\newcommand{\mc}{\mathcal}

\newcommand{\s}{\subseteq}

\newcommand{\vn}{\varnothing}

\renewcommand{\le}{\leqslant}
\renewcommand{\ge}{\geqslant}

\author[I. Colazzo \and E. Jespers \and {\L}. Kubat \and A. Van Antwerpen \and C. Verwimp]
{Ilaria Colazzo \and Eric Jespers \and {\L}ukasz Kubat \and Arne Van Antwerpen \and Charlotte Verwimp}

\address[I. Colazzo (ORCID: 0000-0002-2713-0409)]{Department of Mathematics,
College of Engineering, Mathematics and Physical Sciences, University of Exeter, Exeter EX4 4QF, UK}
\email{I.Colazzo@exeter.ac.uk}

\address[E. Jespers (ORCID: 0000-0002-2695-7949), A. Van Antwerpen (ORCID: 0000-0001-7619-6298),
and C. Verwimp (ORCID: 0000-0003-3128-4854)]{Department of Mathematics and Data Science,
Vrije Universiteit Brussel, Pleinlaan 2, 1050 Brussel}
\email{eric.jespers@vub.be}
\email{arne.van.antwerpen@vub.be}
\email{charlotte.verwimp@vub.be}

\address[{\L}. Kubat (ORCID: 0000-0002-7848-6405)]{University of Warsaw, Institute of Mathematics, Banacha 2, 02-097 Warsaw, Poland}
\email{lukasz.kubat@mimuw.edu.pl}

\title{Finite idempotent set-theoretic solutions of the Yang--Baxter equation}
\subjclass[2010]{Primary: 16T25, 20M25; Secondary: 16W22, 16S36}
\keywords{Yang--Baxter equation, set-theoretic solution, idempotent, structure semigroup, structure algebra.}
\date{}

\begin{document}

\begin{abstract}
It is proven that finite idempotent left non-degenerate set-theoretic solutions $(X,r)$ of the 
Yang--Baxter equation on a set $X$ are
determined by a left simple semigroup structure on $X$ (in particular, a finite union of isomorphic copies of a group)
and some maps $q$ and $\varphi_x$ on $X$, for $x\in X$. This structure turns out to be a group precisely when the associated
structure monoid is cancellative and all the maps $\varphi_x$ are equal to an automorphism of this group. Equivalently,
the structure algebra $K[M(X,r)]$ is right Noetherian, or in characteristic zero it has to be semiprime. The structure
algebra always is a left Noetherian representable algebra of Gelfand--Kirillov dimension one. To prove these results
it is shown that the structure semigroup $S(X,r)$ has a decomposition in finitely many cancellative semigroups $S_u$
indexed by the diagonal, each $S_u$ has a group of quotients $G_u$ that is finite-by-(infinite cyclic) and the union
of these groups carries the structure of a left simple semigroup. The case that $X$ equals the diagonal is fully
described by a single permutation on $X$.
\end{abstract}

\maketitle

\section{Introduction}\label{sec-1}

Let $V$ be a vector space over a field $K$. Recall that a linear map $R\colon V\otimes V\to V\otimes V$ is said to
be a solution of the Yang--Baxter equation if
\[(R\otimes{\id})({\id}\otimes R)(R\otimes{\id})=({\id}\otimes R)(R\otimes{\id})({\id}\otimes R).\]
This equation originates from papers by Baxter \cite{Ba} and Yang \cite{Ya} on statistical physics. The search for
solutions has attracted numerous studies in both mathematical physics and pure mathematics. Topics of interest include
quantum groups, Hopf algebras, Galois theory, group algebras, knot theory, radical rings, pre-Lie algebras, solvable
groups and Garside groups and algebraic structure of quadratic algebras (see for example 
\cite{AngionoVendramin,CJDR,CJORet,Chouraqui,Dehornoy,EPGS,ESS99,FCC,GVdB,KasselBook,LeVe17,RumpRR,RumpClassical,SmokPre,SmokPre1}).
Because of the complexity of the solutions, Drinfeld in \cite{Dri1992} initiated the investigations of set-theoretic
solutions of the Yang--Baxter equation. These are solutions which are induced by a linear extension of a map
$r\colon X\times X\to X\times X$, where $X$ is a basis of $V$. In this case, $r$ satisfies
\[(r\times{\id})({\id}\times r)(r\times{\id})=({\id}\times r)(r\times{\id})({\id}\times r),\]
and one says that $(X,r)$ is a set-theoretic solution of the Yang--Baxter equation (for simplicity, throughout the paper,
we will simply call this a solution). For any $x,y \in X$, we put \[r(x,y)=(\lambda_x(y),\rho_y(x)).\] The solution
is said to be bijective if $r$ is bijective (involutive if $r^{2}=\id$),  left non-degenerate (respectively right non-degenerate)
if all maps $\lambda_x$ (respectively all maps $\rho_x$) are bijective (and non-degenerate if it is left and right non-degenerate).

The work of Etingoff, Soloviev, Schedler \cite{ESS99} and Gateva-Ivanova and Van den Bergh \cite{GVdB} brought renewed
attention to Drinfeld's problem in the particular case of involutive solutions. In \cite{GVdB} the associated structure
monoid $M(X,r)$ and structure algebra $K[M(X,r)]$ (an algebra with defining quadratic relations naturally determined
by the solution $(X,r)$) were investigated, and it was shown that the structure algebra shares many ring and homological
properties of polynomial algebras in commuting variables rovided $X$ also is finite. The ring theoretical properties of
the structure algebra of an arbitrary finite bijective non-degenerate solution have been extensively further investigated,
see for example \cite{GIMa08,MR2927367,MR4278764,MR4213053,MR4024551,MR4105532}. In particular, it is shown that the algebraic
structure determines information of the type of solutions; for example a bijective finite non-degenerate solution is involutive
precisely when its structure algebra is prime (or equivalently, a domain). The introduction of algebraic structures, such
as braces and skew braces \cite{Rump2007, GV17}, breathes new life into the study of involutive and, more generally, bijective
non-degenerate solutions. In fact, any skew brace gives rise to a bijective non-degenerate solution to the Yang--Baxter equation
and conversely, one can always associate a skew brace to such a solution \cite{SmVe18}. In the same vein, in \cite{CoJeVAVe21x},
it is shown that arbitrary left non-degenerate solutions can be studied via the so-called YB-semitrusses (a structure that
specialises to skew braces in the case of bijective solutions).

Jones \cite{MR908150} and Turaev \cite{MR939474} investigated solutions $R\colon V\otimes V \to V\otimes V$ of the Yang--Baxter
equation, where $V$ is a finite-dimensional linear space, which are algebraic of degree two, that is $R^2=aR+b\id$. In case they
are bijective, they give rise to polynomial invariants of knots. In the combinatorial setting, such class of solutions
restricts to set-theoretic solutions $(X,r)$ satisfying $r^2=r$ or $r^2=\id$, since one does not have an addition (as in $\End(V\otimes V)$).
As discussed before, involutive solutions have been widely studied. Whereas, Lebed in \cite{Leb2017,MR4146852}
initiated to investigate idempotent solutions in order to identify the Hochschild cohomology of the plactic monoid and
to show that the quantum symmetrizer yields a quasi-isomorphism between a certain quotient of the braided chain complex for
an idempotent solution, and the Hochschild chain complex for its monoid, with the same coefficients. 

Stanovsk\'{y} and Vojt\v{e}chovsk\'{y} in \cite{StVo21} proved that idempotent left non-degenerate solutions $(X,r)$ are in
one-to-one correspondence with twisted left Ward quasigroups, i.e., left quasigroups $(X,*)$ satisfying the identity $(x*y)*(x*z)=(y*y)*(y*z)$.
The latter are completely described in case $X$ has prime cardinality. They also showed that twisted Ward quasigroups are obtained by
twisting the left division operation in groups (that is $x*y=\psi^{-1}(x^{-1}\cdot y)$ for a group $(X,\cdotp)$), and these
correspond to idempotent latin solutions, i.e., such solutions are determined by a group structure on $X$ and an automorphism of this group.

In this paper, we investigate arbitrary finite left non-degenerate idempotent solutions $(X,r)$ and show that they are determined by a
left simple semigroup structure on $X$, i.e., by finitely many isomorphic finite groups and some mappings on $X$.
Our approach is different and makes use of two results. Firstly, Majid and Gateva-Ivanova \cite{GIMa08} proved that
any solution (without any restriction) can be extended to a solution on its structure monoid.
Hence the latter contains all information on solutions. Secondly, Ced\'o, Jespers and Verwimp \cite{CJV21} showed
that the structure monoid of any left non-degenerate solution is embeddable in the holomorph of the associated
derived structure monoid. The latter is the structure monoid of the associated rack solution. Its special
algebraic structure can then be exploited to obtain information on the original solution. 

We prove that for left non-degenerate finite idempotent solutions $(X,r)$ we have much more rigid
control and, surprisingly, such solutions again are determined by a group and some maps on finitely many isomorphic
copies of this group (Theorem~\ref{solgeneral}). While the structure algebra $K[M(x,r)]$ always is left Noetherian,
it also is right Noetherian only if the number of these isomorphic copies is precisely one, or equivalently the
structure algebra is semiprime in case the field $K$ has characteristic zero (Theorem~\ref{algebrastructure}).

In Section \ref{sec-2}, we recall the necessary background on left non-degenerate solutions $(X,r)$ and the
associated structure monoid $M(X,r)$ and structure semigroup $S=S(X,r)$. In Section \ref{sec-3}, we deal
with such idempotent finite solutions. First, we show that $S$ has a natural decomposition into cancellative
semigroups $S_u$, indexed by elements $u$ in the diagonal $\Lambda =\{\lambda_x^{-1}(x):x\in X\}$. We show
that each $S_u$ is a subsemigroup of its group of quotients $G_u$, which turns out to be a finite conjugacy
group and thus has a finite torsion part, denoted $T(G_u)$. Furthermore, the union of these torsion-groups
is a left simple semigroup and it is in bijective correspondence with $X$. Hence, in the main theorem, we
show that the considered solutions $(X,r)$ are determined by such a left simple semigroup $(X,\cdot)$ and
some maps satisfying some identities. It is then shown that having a unique component $S_u$ (i.e., $S(X,r)$
is embeddable in a group) is equivalent with $(X,\cdot)$ being a group and having an automorphism $\varphi$
so that $\lambda_x (y) =x\cdot \varphi (y)$, and this on its turn is equivalent with its structure algebra
$K[M(X,r)]$ being right Noetherian (or, equivalently, sempiprime if $K$ has characteristic zero).
It turns out that $K[M(X,r)]$ always is a left Noetherian representable algebra (and hence satisfies
a polynomial identity) of Gelfand--Kirillov dimension one. So the algebraic structure of the structure algebra
determines precise information of the solution $r$. We finish the paper with three consequences. First the case
that $|\Lambda |=1$, or equivalently $S$ is embeddable in a group, also corresponds to $(X,r)$ being a latin
solution and such solutions are determined by a group structure on $X$ and an automorphism on this group
(Corollary~\ref{idempotentlatin}). Second the case that $|\Lambda|=|X|$ is simply determined by a permutation
on $X$ (Corollary~\ref{qbijective}). Third, there are only two types of left non-degenerate idempotent
solutions of prime cardinality (Corollary~\ref{primecase}). In Corollary~\ref{idempotentlatin} and
Corollary~\ref{primecase} we recover two results of \cite{StVo21}.

\section{Preliminaries}\label{sec-2}

Throughout we denote by $(X,r)$ a left non-degenerate set-theoretic solution of the Yang--Baxter equation
and we put $r(x,y)=(\lambda_x(y), \rho_y(x))$. So each $\lambda_x \colon X\to X$ is a bijection.

Being a solution of the Yang--Baxter equation translates into the following equalities to hold for all $x,y,z\in X$:
\begin{align}
	\lambda_x(\lambda_y(z)) &=\lambda_{\lambda_x(y)}(\lambda_{\rho_y(x)}(z)),\label{YB1}\\
	\lambda_{\rho_{\lambda_y(z)}(x)}(\rho_z(y)) &=\rho_{\lambda_{\rho_y(x)}(z)}(\lambda_x(y)),\label{YB2}\\
	\rho_z(\rho_y(x)) &=\rho_{\rho_z(y)}(\rho_{\lambda_y(z)}(x)).\label{YB3}
\end{align}
Note that
\[r^2(x,y)=(\lambda_{\lambda_x(y)}(\rho_y(x)),\rho_{\rho_y(x)}(\lambda_x(y))).\]
Hence the solution $(X,r)$ is idempotent, i.e., $r^2=r$, if and only if
\begin{equation}
    \lambda_{\lambda_x(y)}(\rho_y(x))=\lambda_x(y)\qquad\text{and}\qquad\rho_{\rho_y(x)}(\lambda_x(y))=\rho_y(x).\label{r2=r}
\end{equation}

The associated structure monoid of the solution $(X,r)$ is
\[M=M(X,r)=\free{x\in X\mid x \circ y = \lambda_x(y) \circ \rho_y(x),\; x,y\in X}\]
The structure semigroup is \[S=S(X,r)=M(X,r) \setminus \{1\}.\]
The associated derived solution is \[s\colon X^2 \to X^2\colon (x,y) \mapsto (y, \sigma_y(x)),\]
where $\sigma_y(x) = \lambda_y \rho_{\lambda^{-1}_x (y)} (x)$ and its associated structure monoid
(called the derived structure monoid) is denoted $A(X,r)$. Its operation is denoted additively.
Let $\pi\colon M\to A$ be the bijective $1$-cocycle described in \cite{CJV21}. So
\[A=A(X,r)=\free{\pi(x)\text{ for }x\in X\mid \pi(x)+\pi(y)=\pi(y)+\pi(\sigma_y(x)),\; x,y\in X}.\]

Furthermore, \[f\colon M\to A\rtimes{\Img\lambda}\colon m\mapsto (\pi (m), \lambda_{m})\]
is a monoid embedding, where \[\lambda\colon (M,\circ) \to \Aut(A,+)\colon m \mapsto \lambda_m\]
and $\lambda_m =\lambda_x$ if $m=x$. Abusing notation, we will identify $m$ with $f(m)$, i.e.,
we will write $m=(\pi(m),\lambda_m)$. For simplicity reasons, for $a\in A$, we will write $\lambda_a$
for $\lambda_{\pi^{-1}(a)}$. So we simply may write \[M=\{(a,\lambda_a):a\in A\}=\free{x=(x,\lambda_{x}):x\in X}.\]
Therefore, we have a mapping \[\lambda\colon A\to\Aut(A,+)\colon a\mapsto\lambda_a,\]
and, for $a,b\in A$, \[\lambda_{a+\lambda_a(b)}=\lambda_a\circ\lambda_b .\]

We also have a monoid antihomomorphism $\sigma\colon (A,+) \to \End (A,+)\colon a\mapsto \sigma_a$
(``extending'' the maps $\sigma_x$ for $x\in X$). So for $a,b\in A$, we have $a+b=b +\sigma_b (a)$
and thus $A+b\s b+A$. In particular, each right ideal of $A$ is a two-sided ideal of $A$.
For a subset $B$ of $A$ we put, as in \cite{MR4024551}, \[B^e=\{(a,\lambda_a):a\in B\}\s M.\]
Moreover, if $Z$ is a subset of a monoid $T$ then $\free{Z}$ denotes the submonoid of $T$ generated by $Z$.
In case $T$ is a group, $\gr(Z)$ stands for the subgroup of $T$ generated by $Z$.

A fundamental theorem of Gateva-Ivanova and Majid says that set-theoretic solutions are determined
by the solutions associated to their structure monoids.

\begin{thm}[{see \cite[Theorem 3.6]{GIMa08}}]\label{main}
    Assume $(X,r)$ is a solution of the Yang--Baxter equation. If $M=M(X,r)$ then there exist a monoid morphism
    $\lambda\colon(M,\circ)\to\Map(M,M)$ and a monoid antimorphism $\rho\colon(M,\circ)\to\Map(M,M)$, both naturally
    extending the maps $\lambda_x$ and $\rho_x$, respectively, such that, for all $a,b,c\in M$,
    \begin{align}
        \rho_b(c\circ a) &= \rho_{\lambda_a(b)}(c)\circ \rho_b(a), \label{eqrho} \\
        \lambda_b(a \circ c) &= \lambda_b(a) \circ \lambda_{\rho_a(b)}(c),
    \label{eqlambda}
    \end{align} 
    and
    \begin{equation}\label{eqproduct}
        a\circ b=\lambda_{a}(b)\circ\rho_{b}(a).
    \end{equation}
    Let $r_M\colon M\times M\to M\times M$ be defined by $r_M(a,b)=(\lambda_a(b),\rho_b(a))$, for all $a,b\in M$.
    Then $(M,r_M)$ is a set-theoretic solution of the Yang--Baxter equation. Obviously, $r_M$ extends the solution $r$.
\end{thm}

So set-theoretic solutions of the Yang--Baxter equation are determined by the presentation of their structure monoid.
In this paper, we will describe the presentations of the structure semigroup of left non-degenerate idempotent solutions,
and as a consequence, we obtain a description of all such solutions.

So from now on all solutions $(X,r)$ are idempotent and left non-degenerate. The diagonal map
\[q\colon X\to X\colon x\mapsto \lambda_x^{-1}(x)\] will play a crucial role. Put \[\Lambda =\Lambda(X,r)=\Img(q).\]
Clearly $r^2=r$ means that $\rho_y(x) = q(\lambda_{x}(y))$ and thus \[r(x,y)=(\lambda_x(y), q(\lambda_x(y)).\]
Furthermore, the left non-degeneracy yields \[A(X,r)=\free{x\in X\mid x+y =y+y,\;x,y\in X}=\bigcup_{x\in X}\free{x};\]
clearly $\free{x}=\{nx:n\ge 0\}$. Note that if $(X,r)$ is idempotent and left non-degenerate then each element of
$A(X,r)$ is uniquely determined by its length and its last letter; in particular $A(X,r)$ is left cancellative.
We will use this observation freely in the rest of the paper without further reference. The natural length function
on $A$ and $M$ is denoted $|a|$, for $a\in A$, or $a\in M$. All mapping $\lambda_a$ are length preserving.

\section{Main results}\label{sec-3}

In this section, we deal with arbitrary finite left non-degenerate idempotent solutions $(X,r)$. First, we show that
the structure semigroup can be decomposed as a finite union of cancellative semigroups $S_u$, in such a way that it
is graded by $\Lambda$. Next, we prove that each $S_u$ has a group of quotients $G_u$ which is a finite conjugacy group
and thus that its periodic subgroup $T(G_u)$ is finite. This will allow us to obtain a left simple semigroup structure
on $X$ and prove the main results.

\begin{notation}
    Throughout the paper $d$ denotes a multiple of the exponent of the permutation group $\mc{G}(X,r)=\gr(\lambda_x:x\in X)$
    of the solution $(X,r)$. In particular, $\lambda_a^d=\id$ for all $a\in A(X,r)$.
\end{notation}

\begin{prop}\label{decomposition}
    Let $(X,r)$ be a left non-degenerate idempotent solution of the Yang--Baxter equation. Let $S=S(X,r)$
    and $\Lambda=\Lambda(X,r)$. For $u\in\Lambda$ put \[S_u=\{(a,\lambda_a)\in S:\lambda_a^{-1}(a)=|a|u\}.\]
    Then:
    \begin{enumerate}
        \item each $S_u$ is a cancellative subsemigroup of $S$.
        \item $S=\bigcup_{u\in\Lambda} S_u$, a disjoint union.
        \item $S_u \circ S_v \subseteq S_v$ for all $u,v\in\Lambda$.
        \item if $s,t\in S$ then $s\circ t\in S_u$ if and only if $t\in S_u$.
    \end{enumerate}
    In particular, equipping $\Lambda$ with the band structure via $uv=v$ the semigroup $S$
    is a band graded semigroup with homogeneous components that are cancellative semigroups.
    \begin{proof}
        First note that if $u=q(x)\in\Lambda$ then $x\in S_u$. So all $S_u$ are non-empty.
        
        Let $u,v\in\Lambda$ and assume $(a,\lambda_a)\in S_u$ and $(b,\lambda_b)\in S_v$. Then
        \begin{align*}
            \lambda_{a+\lambda_a(b)}^{-1}(a+\lambda_a(b)) & =(\lambda_a\lambda_b)^{-1}(a+\lambda_a(b))
            =\lambda_b^{-1}(\lambda_a^{-1}(a+\lambda_a(b)))\\
            & =\lambda_b^{-1}(\lambda_a^{-1}(a)+b)=\lambda_b^{-1}(\lambda_a^{-1}(a))+\lambda_b^{-1}(b)\\
            & =\lambda_b^{-1}(|a|u)+|b|v=|a+\lambda_a(b)|v. 
        \end{align*}
        Hence $(a,\lambda_a)\circ (b,\lambda_b)=(a+\lambda_a(b), \lambda_{a+\lambda_a(b)})\in S_v$ and thus
        $S_u\circ S_v\s S_v$. In particular, each $S_u$ is a subsemigroup of $S$. Moreover, $S_u\cap S_v=\vn$ for $u\ne v$.
        
        Assume $(a,\lambda_a),(b,\lambda_b),(c,\lambda_c)\in S_u$. If $(a,\lambda_a)\circ (c,\lambda_c)=(b,\lambda_b)\circ(c,\lambda_c)$
        then $a+\lambda_a(c)=b+\lambda_b(c)$ and $\lambda_a=\lambda_b$. This implies $|a|=|b|$ and thus
        \[\lambda_a^{-1}(a)=|a|u=|b|u=\lambda_b^{-1}(b)=\lambda_a^{-1}(b),\] which yields $a=b$. Similarly, we show that 
        $(c,\lambda_c)\circ(a,\lambda_a)=(c,\lambda_c)\circ (b,\lambda_b)$ leads to $a=b$. Hence, $S_u$ is a cancellative semigroup.
        
        To prove that $S=\bigcup_{u\in\Lambda}S_u$, it is enough to show that for each positive integer $n$ and $x\in X$ we have that
        $u=\lambda_{nx}^{-1}(x)\in\Lambda$, because then $\lambda_{nx}^{-1}(nx)=n\lambda_{nx}^{-1}(x)=|nx|u$ and thus
        $(nx,\lambda_{nx})\in S_u$. We prove this by induction on $n$. Clearly, by definition of $\Lambda$, the result holds for $n=1$.
        So, assume that $n>1$ and the result holds for $n-1$. Put $y=\lambda_x^{-1}(x)\in\Lambda\s X$. Because 
        \[\lambda_{nx}=\lambda_{x+\lambda_x(\lambda_x^{-1}((n-1)x))}=\lambda_{x+\lambda_x((n-1)y)}=\lambda_x\lambda_{(n-1)y},\]
        we get \[u=\lambda_{nx}^{-1}(x)=\lambda_{(n-1)y}^{-1}\lambda_x^{-1}(x)=\lambda_{(n-1)y}^{-1}(y)\in\Lambda,\] as desired.
        This finishes the proof of (1), (2) and (3).
        
        Part (4) follows at once from part (3).
    \end{proof}
\end{prop}

\begin{rem}\label{qinlambda}
    In the proof of Proposition~\ref{decomposition} we have shown that $\lambda_{nx}^{-1}(x)\in\Lambda$
    for all $x\in X$ and positive integers $n$.
\end{rem}

\begin{lem}\label{quotientgroup}
    If $u\in\Lambda$ then:
    \begin{enumerate}
        \item $c_u\coloneqq(du,\id) \in S_u$ is a central element in $S_u$.
        \item $S_u =\free{c_u}\circ F_u$, where $F_u\coloneqq\{(a,\lambda_a)\in S_u:0<|a|<d\}$ is a finite set.
        \item $S_u$ has a group of quotients, denoted by $G_u$, that is a finite conjugacy group (FC-group).
        In particular, its elements of finite order form a finite group, denoted by $T(G_u)$.
        \item $G_u/T(G_u)\cong \mb{Z}$ and there exists an element $x_u\in X$ so that $(x_u,\lambda_{x_u})\in S_u$
        and $G_u/T(G_u)$ is generated by $(x_u,\lambda_{x_u}) \circ T(G_u)\in G_u/T(G_u)$.
    \end{enumerate}
    In particular, $G_u \cong T(G_u) \rtimes_{\sigma_u} {\gr(x_u)}\cong T(G_u)\rtimes\mb{Z}$, a semidirect product,
    and the action $\sigma_u$ is given by conjugation by $x_u$.
    
    Furthermore, $S=S(X,r)$ is cancellative if and only if $|\Lambda|=1$.
    \begin{proof}
        As $\lambda_{du}=\id$, we have $\lambda_{du}^{-1}(du)=du$ and thus $c_u\in S_u$. Moreover, if $(nx,\lambda_{nx})\in S_u$
        for some positive integer $n$ and some $x\in X$ then $n\lambda_{nx}(u)=\lambda_{nx} (nu)=nx$. Hence $\lambda_{nx}(u)=x$,
        which leads to
        \begin{align*}
            (nx,\lambda_{nx})\circ c_u & =(nx,\lambda_{nx})\circ(du,\id)=(nx+\lambda_{nx}(du),\lambda_{nx})\\
            & =(nx+d\lambda_{nx}(u),\lambda_{nx})=(nx+dx,\lambda_{nx})=((n+d)x,\lambda_{nx}).
        \end{align*}
        Therefore,
        \begin{align*}
            c_u\circ(nx,\lambda_{nx}) & =(du,\id)\circ (nx,\lambda_{nx})=(du+nx,\lambda_{nx})\\
            & = ((d+n)x,\lambda_{nx})=(nx,\lambda_{nx})\circ c_u.
        \end{align*}
        So part (1) follows.
        
        Let $(nx,\lambda_{nx})\in S_u$ and write $n=kd+l$ for some non-negative integer $k$ and some integer $0<l<d$. Then
        \begin{align*}
            c_u^k\circ(lx,\lambda_{lx}) & =(du,\id)^k\circ(lx,\lambda_{lx})=(kdu,\id)\circ(lx,\lambda_{lx})\\
            & =(kdu+lx,\lambda_{lx})=(kdx+lx,\lambda_{lx})=(nx,\lambda_{lx}).
        \end{align*}
        In particular, $\lambda_{nx}=\lambda_{lx}$. Moreover, $x=\lambda_{nx}(u)=\lambda_{lx}(u)$, which implies
        $(lx,\lambda_{lx})\in F_u$, and part (2) follows.
        
        Next, observe that for any $(a,\lambda_a)\in S_u$, we have
        \begin{equation}
            (a,\lambda_a)^d=(d\lambda_a^{-1}(a),\id)=(d|a|u,\id)=(du,\id)^{|a|}=c_u^{|a|}.\label{2.4-1}
        \end{equation}
        Thus, $G_u\coloneqq S_u\circ\free{c_u}^{-1}$ is the group of quotients of $S_u$ and the central subgroup 
        $\gr(c_u)=\free{c_u}\circ\free{c_u}^{-1}$ of $G_u$ is such that the group $G_u/\gr(c_u)$ is finite (of order
        not exceeding $|F_u|\le d|X|$). Hence $G_u$ has a central (infinite cyclic) subgroup of finite index, and
        thus $G_u$ is an FC-group. Consequently, its torsion subgroup $T(G_u)$ is finite and $G_u/T(G_u)$ is a
        torsion-free abelian group. In particular, part (3) follows.
        
        As $u\in\Lambda$, there exists $x_u\in X$ such that $u=\lambda_{x_u}^{-1}(x_u)$. In particular,
        $x_u=(x_u,\lambda_{x_u})\in S_u$ and $x_u^d=c_u$. Now, suppose $(nx,\lambda_{nx})\in S_u$. Then \eqref{2.4-1}
        implies that $(nx,\lambda_{nx})^d=c_u^n=x_u^{dn}=(x_u^n)^d$. As the group $G_u/T(G_u)$ is torsion-free abelian,
        it follows that $(nx,\lambda_{nx})\in T(G_u)\circ x_u^n$. So, indeed
        $G_u\cong T(G_u)\rtimes_{\sigma_u}{\gr((x_u,\lambda_{x_u}))}\cong T(G_u)\rtimes\mb{Z}$ and the action of
        $\sigma_u$ is conjugation by $(x_u,\lambda_{x_u})$. This finishes the proof of part (4).
        
        The last part of the result follows from the fact that if $(x,\lambda_x)\in S_u$ and $(y,\lambda_y)\in S_v$,
        with $x,y\in X$ and distinct $u,v\in\Lambda$ then
        $(x,\lambda_x)^{d} \circ (y,\lambda_y)^{d} =c_u \circ c_v =c_v^{2}=(y,\lambda_y)^d\circ (y,\lambda_y)^{d}$,
        while $(x,\lambda_x)^{d}=c_u \ne (y,\lambda_y)^d=c_v$.
    \end{proof}
\end{lem}

Next, we describe the periodic elements of $G_u$. To do so we need to introduce some notation, which will be useful later as well.
Namely, for $u\in\Lambda$ define \[X_u\coloneqq\{x\in X:\lambda_{dx}(u)=x\}=\{x\in X:(dx,\lambda_{dx})\in S_u\}.\]
As $S(X,r)=\bigcup_{u\in\Lambda}S_u$, it is clear that we have a decomposition $X=\bigcup_{u\in\Lambda}X_u$.
Moreover, for $x\in X_u$ put \[t_{u,x}\coloneqq(dx,\lambda_{dx})\circ(du,\id)^{-1}=(dx,\lambda_{dx})\circ c_u^{-1}\in G_u.\]

\begin{lem}\label{lemmatorsion}
    If $u\in\Lambda$ then $T(G_u)=\{t_{u,x}:x\in X_u\}$. Furthermore, if $x,y\in X_u$ then:
    \begin{enumerate}
        \item $t_{u,u}$ is the identity of $T(G_u)$.
        \item $t_{u,x}\ne t_{u,y}$ provided $x\ne y$. In particular, $|T(G_u)|=|X_u|$.
        \item $t_{u,x}\circ t_{u,y}=t_{u,\lambda_{dx}(y)}$.
        \item $t_{u,x}=(kdx,\lambda_{kdx})\circ(kdu,\id)^{-1}$ for each positive integer $k$.
        \item the order of the element $t_{u,x}$ is a divisor of $d$.
        \item $\lambda_x=\lambda_{dx}\lambda_u$ and if $q(u)=u$ then $t_{u,x}=(x,\lambda_x) \circ (u,\lambda_u)^{-1}$.
        \item $X=\bigcup_{u\in\Lambda}(\{ (x,\lambda_x): x\in X \} \cap S_u)$  and if $(x,\lambda_x)\in S_u$
        then $(x,\lambda_x)\circ(x_u,\lambda_{x_u})^{-1}\in T(G_u)$.
    \end{enumerate}
    \begin{proof}
        Because $\lambda_{du}=\id$, we get $u\in X_u$ and $t_{u,u}=1$. So part (1) is obvious. Since $dx\ne dy$, part (2) also follows.
        
        Next, write as before $c_u=(du,\id)\in G_u$ and observe that
        \begin{align*}
            t_{u,x} \circ t_{u,y} &=(dx,\lambda_{dx})\circ c_u^{-1} \circ (dy,\lambda_{dy})\circ c_u^{-1}
            =(dx+d\lambda_{dx}(y),\lambda_{dx}\lambda_{dy})\circ c_u^{-2}\\
            & =(du+d\lambda_{dx}(y),\lambda_{dx}\lambda_{dy})\circ c_u^{-2}
            =c_u\circ(d\lambda_{dx}(y),\lambda_{dx}\lambda_{dy})\circ c_u^{-2}\\
            & =(d\lambda_{dx}(y),\lambda_{d\lambda_{dx}(y)})\circ c_u^{-1}=t_{u,\lambda_{dx}(y)}.
        \end{align*}
        Hence part (3) follows as well.
        
        Choose $(a,\lambda_a)\in S_u$ with $a=kdx$ for some positive integer $k$ and some $x\in X$. Then
        \begin{align}
            \begin{aligned}
                c_u\circ(a,\lambda_a) & =(du,\id)\circ(kdx,\lambda_{kdx})=(du+kdx,\lambda_{kdx})=((k+1)dx,\lambda_{kdx}),\\
                c_u^k\circ(dx,\lambda_{dx}) & =(kdu,\id)\circ(dx,\lambda_{dx})=(kdu+dx,\lambda_{dx})=((k+1)dx,\lambda_{dx}).
            \end{aligned}\label{2.5-1}
        \end{align}
        In particular, $\lambda_{dx}=\lambda_{kdx}=\lambda_a$. Hence $\lambda_{dx}(u)=\lambda_a(u)=u$ because $(a,\lambda_a)\in S_u$,
        and thus $x\in X_u$. It also follows from \eqref{2.5-1} that $(a,\lambda_a)\circ c_u^{-k}=(dx,\lambda_{dx})\circ c_u^{-1}=t_{u,x}$.
        So part (4) follows. As $S_u$, and in consequence, its central localization $G_u=S_u\free{c_u}^{-1}$, inherits gradation by length
        from $S$, we know that elements of $T(G_u)$ are of length zero, and thus each such element may be written in the form
        $(a,\lambda_a)\circ c_u^{-k}$, where $(a,\lambda_a)\in S_u$ is of length $kd$. This yields $T(G_u)=\{t_{u,x}:x\in X_u\}$.
        
        Since $(dx,\lambda_{dx})^d=c_u^d$ (see \eqref{2.4-1} in the proof of Lemma \ref{quotientgroup}), we get $t_{u,x}^d=1$.
        In particular, part (5) follows.
        
        We prove part (6). Because $\lambda_{dx}(u)=x$, we obtain $du+x=dx+\lambda_{dx}(u)$, and thus
        \begin{align*}
            c_u\circ(x,\lambda_x) & =(du,\id)\circ(x,\lambda_x)=(du+x,\lambda_x)=(dx+\lambda_{dx}(u),\lambda_x)\\
            & =(dx,\lambda_{dx}) \circ (u,\lambda_u)=(dx+\lambda_{dx}(u),\lambda_{dx}\lambda_u).
        \end{align*}
        Therefore, $\lambda_x=\lambda_{dx}\lambda_u$. Furthermore, if $q(u)=u$ then $(u,\lambda_u)\in S_u$ and thus
        $t_{u,x}=(dx,\lambda_{dx})\circ c_u^{-1}=(x,\lambda_{x}) \circ (u,\lambda_u)^{-1}$. Hence (6) follows.
        
        Finally part (7) follows from the previous grading argument, which says that $T(G_u)$ consists of the elements of degree $0$.
    \end{proof}
\end{lem}

\begin{rem}\label{fixedpoint}
    We know that if $(X,r)$ is a left non-degenerate idempotent solution of the Yang--Baxter equation and $a,b\in A(X,r)$
    are of equal length then $\lambda_a=\lambda_b$ or $\lambda_a\lambda_b^{-1}$ has no fixed-points. Indeed, if
    $\lambda_a(x)=\lambda_b(x)$ for some $x\in X$ then
    \[\lambda_a\lambda_x=\lambda_{a+\lambda_a(x)}=\lambda_{b+\lambda_b(x)}=\lambda_b\lambda_x,\]
    and thus $\lambda_a=\lambda_b$. Since $\lambda_{dx}=\lambda_x\lambda_u^{-1}$ for $u\in\Lambda$ and $x\in X_u$,
    it follows $\lambda_{dx}=\id$ or $\lambda_{dx}$ does not have fixed-points.
\end{rem}

Recall that a finite simple semigroup is isomorphic with a matrix semigroup $\mc{M}(G,n,m,P)$ consisting of elements
$(g,i,j)$ with $g\in G$, $1\le i\le n$ and $1\le j\le m$. The multiplication is given by the $m\times n$ sandwich matrix
$P=(p_{ji})$ with entries in $G$ as follows: $(g,i,j)(h,k,l) =(gp_{jk}h,i,l)$. Conversely, any such semigroup (even if $G$
is infinite) is a simple semigroup. Furthermore, such a semigroup is left cancellative if and only if $n=1$ and in such a
case, one may assume that each entry of $P$ is the identity element of the group. We will denote such a matrix as $J$.
Clearly $\mc{M}(G,1,n,J) =\bigcup_{j=1}^{m}G_j$, where $G_j=\{(g,1,j):g\in G\}$. Moreover, in this case, $G_iG_j=G_j$ and
all $G_i$ are isomorphic and they are the maximal subgroups. Note that every idempotent of this semigroup is a left identity.
For details we refer the reader to \cite{CP}. 

We now prove that \[T(G^e)\coloneqq\bigcup_{u\in\Lambda}T(G_u)\]
is equipped with a semigroup structure that extends each group $T(G_u)$.

\begin{lem}\label{infotorsioncover}
    $T(G^e)$ is a left cancellative simple semigroup for the binary operation $\circ$ defined as
    \begin{equation}
        t_{u,x}\circ t_{v,y}\coloneqq t_{v,\lambda_{dx}(y)}\label{defproducttorsion}
    \end{equation}
    for $u,v\in\Lambda$, $x\in X_u$ and $y\in X_v$. Its maximal subgroups are the groups $T(G_u)$. Hence:
    \begin{enumerate}
        \item $T(G_u) \cong T(G_v)$ for all $u,v\in\Lambda$.
        \item $T(G^e) \cong \mc{M}(T(G_u),1,|\Lambda|,J)$ for some (and thus any) $u\in\Lambda$.
        \item $|X| = |\Lambda|\cdot |T(G_u)|$ for any $u\in\Lambda$.
        \item $X=\bigcup_{u\in\Lambda}X_u$ is a disjoint union and $|X_u|=|X_v|$ for all $u,v\in\Lambda$.
    \end{enumerate}
    \begin{proof}
        First notice that
        \[(dx,\lambda_{dx})\circ(dy,\lambda_{dy})=(dx+\lambda_{dx}(dy),\lambda_{dx}\lambda_{dy})=(du+\lambda_{dx}(dy),\lambda_{dx}\lambda_{dy})\]
        and also that \[(du,\id)\circ(\lambda_{dx}(dy),\lambda_{\lambda_{dx}(dy)})=(du+\lambda_{dx}(dy),\lambda_{\lambda_{dx}(dy)}).\]
        Hence, it follows that
        \begin{equation}
            \lambda_{dx} \lambda_{dy} =\lambda_{\lambda_{dx}(dy)}=\lambda_{d\lambda_{dx}(y)}.\label{2.7-1}
        \end{equation}
        In particular, $\lambda_{d\lambda_{dx}(y)}(v)=\lambda_{dx}\lambda_{dy}(v)=\lambda_{dx}(y)$,
        which assures that $\lambda_{dx}(y)\in X_v$ and thus $t_{v,\lambda_{dx}(y)}\in G_v$.
        
        Next we verify that the operation \eqref{defproducttorsion} is associative. Hence, if additionally
        $w\in\Lambda$ and $z\in X_w$ then using \eqref{2.7-1} we get
        \begin{align*}
            t_{u,x}\circ(t_{v,y}\circ t_{w,z}) & =t_{u,x}\circ t_{w,\lambda_{dy}(z)}=t_{w,\lambda_{dx}\lambda_{dy}(z)}\\
            & =t_{w,\lambda_{d\lambda_{dx}(y)}(z)}=t_{v,\lambda_{dx}(y)}\circ t_{w,z}=(t_{u,x}\circ t_{v,y})\circ t_{w,z}.
        \end{align*}
        Consequently, we have shown that $T(G^e)$ is semigroup that contains each $T(G_u)$ as a subgroup and $T(G_u)\circ T(G_v)=T(G_v)$
        for all $u,v\in\Lambda$. It follows that $T(G^e)$ is a left cancellative simple semigroup with maximal subgroups $T(G_u)$.
        Hence all $T(G_u)$ are isomorphic, and in particular they all have the same cardinality and
        \[T(G^e)\cong\mc{M}(H,1,|\Lambda|,J),\] where $H$ is a finite group isomorphic to each $T(G_u)$.
        The remaining statements of the result are now obvious.
    \end{proof}
\end{lem}

The next step is to prove that \[G^e\coloneqq\bigcup_{u\in\Lambda}G_u\] may be equipped with a semigroup
structure so that each $G_u$ is a subgroup and $S(X,r)$ is a subsemigroup. From Lemma~\ref{quotientgroup}
we know that an arbitrary element of $G_u$ can be written in the form
$(a,\lambda_a)\circ c_u^n$ with $(a,\lambda_a)\in S_u$ and $n\in\mb{Z}$. Recall that $c_u=(du,\id)\in Z(S_u)$.

\begin{thm}\label{descriptionsimplecover}
    $G^e$ is left cancellative simple semigroup for the binary operation $\circ$ defined as
    \[((a,\lambda_a)\circ c_u^n)\circ((b,\lambda_b)\circ c_v^m)=(a+\lambda_a(b),\lambda_{a+\lambda_a(b)})\circ c_v^{n+m}\]
    for $u,v\in\Lambda$, $(a,\lambda_a)\in S_u$, $(b,\lambda_b)\in S_v$ and $n,m\in\mb{Z}$. Moreover, $S=S(X,r)$ is the
    subsemigroup of
    \[G^e=\bigcup_{u\in \Lambda}G_u =\bigcup_{u\in \Lambda}(T(G_u)\rtimes_{\sigma_u}\gr((x_u,\lambda_{x_u})))
    =\bigcup_{u\in \Lambda}(T(G_u)\rtimes_{\sigma_{u}}\mb{Z})\cong\mc{M}(G_u,1,|\Lambda|,J),\]
    generated by the elements $t\circ (x_{u},\lambda_{x_u})$, with $t\in T(G_u)$ and $u\in \Lambda$.
    Furthermore, \[T(G_u)=\{(x,\lambda_x)\circ(x_u,\lambda_{x_u})^{-1}:x\in X \text{ with } q(x)=u\}.\]
    In particular, the structural algebra $K[S(X,r)]$ is a subalgebra of the $\Lambda$-graded ring
    \[\bigoplus_{u\in \Lambda} K[G_{u}]=\bigoplus_{u\in \Lambda}K[T(G_u)][t_u^{\pm 1};\sigma_u],\]
    where each summand is a skew Laurent polynomial algebra over the group algebra $K[T(G_u)]$ of the finite group $T(G_u)$.
    \begin{proof}
        Clearly, if the proposed binary operation yields a semigroup then it contains $S$ as a subsemigroup
        (as one then multiplies elements with $n=m=0$ precisely as in $S$).
        
        We need to prove that the operation is well-defined. Notice that 
        \[(a+\lambda_a(b),\lambda_{a+\lambda_a(b)})=(a,\lambda_a)\circ(b,\lambda_b)\in S_u\circ S_v\s S_v,\]
        and thus indeed $(a+\lambda_a(b),\lambda_{a+\lambda_a(b)})\circ c_v^{n+m}\in G_v$.
        
        Now, we verify that the product $\circ$ defined in $G^e$ is independent of representatives in $G_u$ and $G_v$.
        So, suppose first that $(a,\lambda_a)\circ c_u^n=(a',\lambda_{a'})\circ c_u^{n'}$ for some $(a',\lambda_{a'})\in S_u$
        and some $n'\in\mb{Z}$. Without loss of generality we may assume that $n\ge n'$. Then, in $S_u$, we have
        $(a',\lambda_{a'})=(a,\lambda_a)\circ c_u^{n-n'}=c_u^{n-n'}\circ (a,\lambda_a)$. In particular,
        $\lambda_{a'}=\lambda_a$ and $a'=(n-n')du+a$. Since $(a+\lambda_a(b),*)\in S_v$, we get
        \begin{align*}
            (a'+\lambda_{a'}(b),*)\circ c_v^{n'+m} & =((n-n')du+a+\lambda_a(b),*)\circ c_v^{n'+m}\\
            & =c_u^{n-n'}\circ(a+\lambda_a(b),*)\circ c_v^{n'+m}\\
            & =c_v^{n-n'}\circ(a+\lambda_a(b),*)\circ c_v^{n'+m}\\
            & =(a+\lambda_a(b),*)\circ c_v^{n-n'}\circ c_v^{n'+m}\\
            & =(a+\lambda_a(b),*)\circ c_v^{n+m}.
        \end{align*}
        So the product is independent of the choice of the representative of the left factor.
        
        Next, we verify the independence of the right factor. So, suppose $(b,\lambda_b)\circ c_v^m=(b',\lambda_{b'})\circ c_v^{m'}$
        for some $(b',\lambda_{b'})\in S_v$ and some $m'\in\mb{Z}$. Again, without loss of generality, we may assume that $m\ge m'$.
        Then, in $S_v$, we have $(b',\lambda_{b'})=(b,\lambda_b)\circ c_v^{m-m'}=c_v^{m-m'}\circ(b,\lambda_b)$. In particular,
        $b'=(m-m')dv+b$. Since $(a+\lambda_a(b),*)\in S_v$, we get
        \begin{align*}
            (a+\lambda_a(b'),*)\circ c_v^{n+m'} & =(a+\lambda_a((m-n')dv+b),*)\circ c_v^{n+m'}\\
            & =(a+(m-m')d\lambda_a(v)+\lambda_a(b),*)\circ c_v^{n+m'}\\
            & =((m-m')dv+a+\lambda_a(b),*)\circ c_v^{n+m'}\\
            & =c_v^{m-m'}\circ(a+\lambda_a(b),*)\circ c_v^{n+m'}\\
            & =(a+\lambda_a(b),*)\circ c_v^{m-m'}\circ c_v^{n+m'}\\
            & =(a+\lambda_a(b),*)\circ c_v^{n+m}.
        \end{align*}
        So the product is indeed well-defined. Notice that the product of $(a,\lambda_a)\circ c_u^n$
        and $(b,\lambda_b)\circ c_v^m$ is determined by the product of $(a,\lambda_a)$ and $(b,\lambda_b)$
        in $S$ and the product of the elements $c_u^n\circ c_v^m=c_v^{n+m}$. Clearly this definition of the
        product on the elements of type $c_u^n$, for various $u\in\Lambda$ and $n\in\mb{Z}$, is associative.
        Hence, the proposed product on $G^e$ yields an associative product.
        
        It is readily verified that the elements $t_{u,u}$ are left identities of $G^e$. Since each $G_u$ is a
        group we obtain that $G^e$ is a left cancellative simple semigroup with maximal subgroups the groups $G_u$;
        in particular all groups $G_u$ are isomorphic. So $G^e \cong\mc{M}(H,1,|\Lambda|,J)$, where $H$ is a group
        that is isomorphic to $G_u\cong T(G_u)\rtimes{\gr(x_u)}$ and $\gr(x_u)\cong\mb{Z}$. It follows that
        \[G^e \cong \mc{M}(G_u,1,|\Lambda|,J).\] Next we prove that the generators $(x,\lambda_x)$ are of the
        form $t\circ (x_u,\lambda_{x_u})$, with $t\in G_u$. Indeed, let $u\in \Lambda$ be such that
        $(x,\lambda_x)\in S_u$, i.e., $q(x)=u$. By Lemma~\ref{lemmatorsion},
        $t\coloneqq(x,\lambda_x)\circ(x_u,\lambda_{x_u})^{-1}\in T(G_u)$.
        Hence $(x,\lambda_x) =t\circ (x_{u},\lambda_{x_u})\in T(G_u)\circ (x_u,\lambda_{x_u})$. Finally,
        since the cardinality of $\bigcup_{u\in \Lambda} T(G_u)\circ (x_u, \lambda_{x_u})$ is the same as
        the cardinality of $X$, it follows that all elements in the union belong to the generating set
        $\{(x,\lambda_x):x\in X\}$ of $S(X,r)$ and the result follows.
    \end{proof}
\end{thm}

Recall that $r(x,y) =(\lambda_x (y), q(\lambda_x(y))$, for $x,y\in X$. From Lemma~\ref{lemmatorsion}
we also know that $\lambda_x = \lambda_{dx}\lambda_u$ for $x\in X_u$. Furthermore, from Lemma~\ref{infotorsioncover},
$t_{u,x} \circ t_{v,y}=t_{v,\lambda_{dx}(y)}$. Hence, via the bijective mapping (see Lemma~\ref{infotorsioncover}) 
\[f\colon X\to T(G^e)\colon x\mapsto t_{u,x}\text{ if $x\in X_u$ for some }u\in\Lambda,\]
we obtain a semigroup structure $(X,\cdotp)$ on $X$ for the binary operation $\cdot$ defined for $x,y\in X$ as
\[x\cdot y =\lambda_{dx}(y).\] Then $(X,\cdotp)$ is a simple left cancellative semigroup with maximal subgroups
isomorphic with $f^{-1}(G_u)$ for $u\in\Lambda$. Note that the elements $u\in\Lambda$ correspond to the idempotents
$t_{u,u}$ of $T(G^e)$ and hence the elements $u\in\Lambda$ are left identities in $(X,\cdotp)$. Furthermore, for
$x\in X_u$ and $y\in X$, \[\lambda_x(y) = x\cdot \lambda_u (y).\] Clearly we then get, for $x\in X_u$ and $y\in X$,
\[r(x,y) =(x\cdot \lambda_u (y),q(x\cdot \lambda_u (y)) ).\]

\begin{lem}\label{qq}
    We have $q^k(x)=\lambda_{kx}^{-1}(x)$ for each $x\in X$ and each positive integer $k$. Moreover, $q^{d+1}=q$.
    \begin{proof}
        We prove this by induction on $k$. For $k=1$ this is clear. Assume now that it holds for $k$.
        Put $y=\lambda_{kx}^{-1}(x)=q^k(x)\in X$. Then we have $\lambda_{(k+1)x}=\lambda_{kx+\lambda_{kx}(y)}=\lambda_{kx}\lambda_y$
        and thus \[\lambda_{(k+1)x}^{-1}(x)=\lambda_y^{-1}\lambda_{kx}^{-1}(x)=\lambda_y^{-1}(y)=q(y)=q^{k+1}(x),\] as desired.

        Next, notice that if $u=q^d(x)\in X$ then we know by the first paragraph that $u\in\Lambda$ and $x\in X_u$.
        Thus, by Lemma~\ref{lemmatorsion}, we get $\lambda_x=\lambda_{dx}\lambda_u$.
        Therefore, $\lambda_{(d+1)x}=\lambda_{dx+\lambda_{dx}(u)}=\lambda_{dx}\lambda_u=\lambda_x$. Consequently,
        $q^{d+1}(x)=\lambda_{(d+1)x}^{-1}(x)=\lambda_x^{-1}(x)=q(x)$.
    \end{proof}
\end{lem}

Note that Lemma \ref{qq} implies that the map $q$ is invertible if and only if $q^d=\id$, or equivalently $\Lambda=X$.

We introduce one more notation. Namely, for $x\in X$, put \[\varphi_x\coloneqq\lambda_{q^d(x)}.\]
In particular, $\varphi_x=\varphi_u$ for each $x\in X_u$. So \[\lambda_x (y) = x\cdot\lambda_{q^d(x)}(y)=x\cdot \varphi_x(y).\]
Hence the solution $r$ has the form \[r(x,y)=(x\cdot\varphi_x(y),q(x\cdot\varphi_x(y))).\]

As $r^2=r$, the idempotent equalities \eqref{r2=r} become
\[\rho_y(x)=q(\lambda_x(y))\qquad\text{and}\qquad q(\lambda_{\lambda_x(y)}(q(\lambda_x(y))))=q(\lambda_x(y)).\]
As $z=\lambda_x(y)$ may be an arbitrary element of $X$, the latter equality reduces to \[q(\lambda_z(q(z)))=q(z).\]

Note that 
\begin{align*}
    \lambda_x\lambda_y(z) & =x\cdot\varphi_x(y\cdot\varphi_y(z)),\\
    \lambda_{\lambda_x(y)}\lambda_{\rho_y(x)}(z) & =\lambda_{\lambda_x(y)}\lambda_{q(\lambda_x(y))}(z)\\
    & =x\cdot\varphi_x(y)\cdot\varphi_{x\cdot\varphi_x(y)}(q(x\cdot\varphi_x(y))\cdot\varphi_{q(x\cdot\varphi_x(y))}(z)).
\end{align*}
Since the semigroup $(X,\cdotp)$ is left cancellative and because elements of $\Lambda$ are left identities for $(X,\cdotp)$,
the Yang--Baxter equation \eqref{YB1} becomes
\[\varphi_x(y\cdot\varphi_y(z))=\varphi_x(y)\cdot\varphi_{x\cdot\varphi_x(y)}\varphi_{q(x\cdot\varphi_x(y))}(z).\]
Note that if all $\varphi_x =\varphi$ (for example if $\Lambda$ is a singleton then this equation becomes
$\varphi(y\cdot\varphi(z))=\varphi(y)\cdot\varphi(\varphi(z))$ and thus writing $t=\varphi(z)$ we obtain
$\varphi(y\cdot t)=\varphi(y)\cdot\varphi(t)$, i.e., $\varphi$ is a homomorphism of the semigroup $(X,\cdotp)$.

Similarly, taking into account that
\begin{align*}
    \lambda_{\rho_{\lambda_y(z)}(x)}(\rho_z(y))
    & =\lambda_{q(\lambda_x\lambda_y(z))}(q(\lambda_y(z)))\\
    & =q(x\cdot\varphi_x(y\cdot\varphi_y(z)))\cdot\varphi_{q(x\cdot\varphi_x(y\cdot\varphi_y(z)))}(q(y\cdot\varphi_y(z))),\\
    \rho_{\lambda_{\rho_y(x)}(z)}(\lambda_x(y)) & =q(\lambda_{\lambda_x(y)}\lambda_{\rho_y(x)}(z))\\
    & =q(x\cdot\varphi_x(y)\cdot\varphi_{x\cdot\varphi_x(y)}(q(x\cdot\varphi_x(y))\cdot\varphi_{q(x\cdot\varphi_x(y))}(z))),
\end{align*}
one can check that the Yang--Baxter equation \eqref{YB2} simplifies to
\[\varphi_{q(x\cdot\varphi_x(y\cdot\varphi_y(z)))}(q(y\cdot\varphi_y(z)))=
q(x\cdot\varphi_x(y)\cdot\varphi_{x\cdot\varphi_x(y)}(q(x\cdot\varphi_x(y))\cdot\varphi_{q(x\cdot\varphi_x(y))}(z))).\]
The latter simplifies further because elements $\Lambda$ are left identities to
\[\varphi_{q(x\cdot\varphi_x(y\cdot\varphi_y(z)))}(q(y\cdot\varphi_y(z)))=
q(x\cdot\varphi_x(y)\cdot\varphi_{x\cdot\varphi_x(y)}(\varphi_{q(x\cdot\varphi_x(y))}(z))).\]

In case $\Lambda$ is a singleton then this equation is trivially satisfied, and it thus not give any extra information.

Finally, since
\begin{align*}
    \rho_z(\rho_y(x))
    & =q(\lambda_{\rho_y(x)}(z))\\
    & =q(\lambda_{q(\lambda_x(y))}(z))\\
    & =q(q(x\cdot\varphi_x(y))\cdot\varphi_{q(x\cdot \varphi_x (y))}(z)),\\
    \rho_{\rho_z(y)}(\rho_{\lambda_y(z)}(x))
    & =q(\lambda_{\rho_{\lambda_y(z)}(x)}(\rho_z(y)))\\
    & =q(\lambda_{q(\lambda_x\lambda_y(z))}(q(\lambda_y(z))))\\
    & =q(q(\lambda_x\lambda_y(z))\cdot\varphi_{q(\lambda_x\lambda_y(z))}(q(\lambda_y(z))))\\
    & =q(q(x\cdot\varphi_x(y\cdot\varphi_y(z)))\cdot\varphi_{x\cdot\varphi_x(y\cdot\varphi_y(z))}(q(y\cdot\varphi_y(z)))),
\end{align*}
the Yang--Baxter equation \eqref{YB3} is equivalent to
\[q(\varphi_{q(x\cdot\varphi_x(y))}(z))= q(\varphi_{x\cdot\varphi_x(y\cdot\varphi_y(z))}(q(y\cdot\varphi_y(z)))).\]
Also this equation is trivially satisfied if $\Lambda$ is a singleton.

Hence we have shown the following theorem.

\begin{thm}\label{solgeneral}
    Assume $(X,r)$ is a finite left non-degenerate idempotent solution of the Yang--Baxter equation.
    Then the following properties are satisfied:
    \begin{enumerate}
        \item there exists a binary operation $\cdot$ on $X$ so that $(X,\cdotp)$ is a left cancellative simple semigroup.
        \item there exits a surjective mapping $q\colon X\to\Lambda$, where $\Lambda$ is the set of idempotents of $(X,\cdotp)$.
        \item for each $x\in X$ there exists a permutation $\varphi_x\in\Sym(X)$ so that the following identities hold:
        \begin{align}
            \varphi_x(y\cdot\varphi_y(z)) & =\varphi_x(y)\cdot\varphi_{x\cdot\varphi_x(y)}\varphi_{q(x\cdot\varphi_x(y))}(z),\label{fineq1}\\
            \varphi_{q(x\cdot\varphi_x(y\cdot\varphi_y(z)))}(q(y\cdot\varphi_y(z)))
            & =q(x\cdot\varphi_x(y)\cdot\varphi_{x\cdot\varphi_x(y)}(\varphi_{q(x\cdot\varphi_x(y))}(z))),\label{fineq2}\\
            q(\varphi_{q(x\cdot\varphi_x(y))}(z)) & =q(\varphi_{x\cdot\varphi_x(y\cdot\varphi_y(z))}(q(y\cdot\varphi_y(z)))),\label{fineq3}\\
            q(x\cdot\varphi_x(q(x))) & =q(x).\label{fineq4}
        \end{align}
        \item the solution $r$ has the form
        \begin{equation}
            r(x,y)=(x\cdot\varphi_x(y),q(x\cdot\varphi_x(y))).\label{sol-form}
        \end{equation}
    \end{enumerate}
    Conversely, if $(X,\cdotp)$ is a finite left cancellative simple semigroup with the set of idempotents $\Lambda$,
    and there exist a surjective mapping $q\colon X\to\Lambda$ and permutations $\varphi_x\in\Sym(X)$, for $x\in X$,
    satisfying conditions \eqref{fineq1}--\eqref{fineq4}, then the formula \eqref{sol-form} yields a finite left
    non-degenerate idempotent solution $(X,r)$ of the Yang--Baxter equation.
    \begin{proof}
        The necessity has been proven before the statement. The proof of sufficiency can be obtained by reversing
        the order in argumentation. Namely, define $\lambda_x(x)=x\cdot\varphi_x(y)$ and $\rho_y(x)=q(x\cdot\varphi_x(y))$
        for $x,y\in X$. Because the maps $q$ and $\varphi_x$ for $x\in X$ satisfy \eqref{fineq1}, \eqref{fineq2} and
        \eqref{fineq3}, it follows that the equalities \eqref{YB1}, \eqref{YB2} and \eqref{YB3} hold. Hence $(X,r)$
        is a solution of the Yang--Baxter equation of the form \eqref{sol-form}. One also easily verifies that the equality
        \eqref{fineq4} translates to $r^2=r$. Finally, as the semigroup $X$ is left cancellative and $\varphi_x$ is bijective,
        the map $\lambda_x$ is injective, and thus bijective as $X$ is finite. Therefore, the proof is complete.
    \end{proof}
\end{thm}

In the next main result we show that $\Lambda$ is a singleton, i.e., $S(X,r)$ is contained in a group, precisely when
the structure algebra $K[M(X,r)$ is right Noetherian. So also for idempotent left non-degenerate solutions the structure
algebra determines crucial information of solutions. First we show that structure algebras of finite left non-degenerate
idempotent solutions are always left Noetherian.

\begin{prop}\label{stralg}
    Assume $(X,r)$ is a finite left non-degenerate idempotent solution of the Yang--Baxter equation. If $M=M(X,r)$ and $K$
    is a field then $K[M]$ is a finite module over a $K$-subalgebra that is isomorphic to a polynomial algebra in one variable.
    In particular, $K[M]$ is a left Noetherian representable algebra (and hence satisfies a polynomial identity)
    of Gelfand--Kirillov dimension one.
    \begin{proof}
        Fix, arbitrarily, an element $t\in X$. We claim that each element $w\in S(X,r)$ may be written (uniquely) as 
        $w=(t,\lambda_t)^{n-1}\circ(x,\lambda_x)$ for some $x\in X$ and some positive integer $n$. We shall prove
        this by induction on $n=|w|$, the length of $w$. If $n=1$ then this is obvious. So, let $n>1$
        and assume the result holds for elements of $S(X,r)$ of length $<n$. Write $w=v\circ(z,\lambda_z)$
        for some $v\in S(X,r)$ of length $n-1$ and some $z\in X$. Then, by the induction hypothesis,
        $v=(t,\lambda_t)^{n-2}\circ(y,\lambda_y)$ for some $y\in X$. Let $x=\lambda_t^{-1}(\lambda_y(z))\in X$.
        Then, taking into account that $t+\lambda_y(z)=y+\lambda_y(z)$ in $A(X,r)$, we obtain
        \begin{align*}
            (t,\lambda_t)\circ(x,\lambda_x) & =(t+\lambda_t(x),\lambda_{t+\lambda_t(x)})=(t+\lambda_y(z),\lambda_{t+\lambda_y(z)})\\
            & =(y+\lambda_y(z),\lambda_{y+\lambda_y(z)})=(y,\lambda_y)\circ(z,\lambda_z).
        \end{align*}
        Therefore,
        \begin{align*}
            w & =v\circ(z,\lambda_z)=(t,\lambda_t)^{n-2}\circ(y,\lambda_y)\circ(z,\lambda_z)\\
            & =(t,\lambda_t)^{n-2}\circ(t,\lambda_t)\circ(x,\lambda_x)=(t,\lambda_t)^{n-1}\circ(x,\lambda_x),
        \end{align*}
        as claimed. In particular, \[K[M]=R+\sum_{x\in X}R\circ(x,\lambda_x)\] is a finite left $R$-module over the polynomial
        algebra in one variable $R=K[(t,\lambda_t)]$. Thus, clearly, $K[M]$ is a left Noetherian algebra of Gelfand--Kirillov
        dimension equal to one (see, e.g., \cite{MR}), while the representability of $K[M]$ is a consequence of \cite{Ana1989}.
    \end{proof}
\end{prop}

{It is worth mentioning that as a consequence of the previous result and a  result of Anan'in \cite{Ana1989} the structure
monoid $M(X,r)$ of a finite left non-degenerate idempotent solution $(X,r)$ is a linear monoid, i.e., a submonoid of the
multiplicative semigroup of a matrix ring over a field, and thus
\cite[Proposition 1, p.~221, Proposition 7, p.~280--281, and Theorem 14, p.~284]{Okn1991} yields
\[\clK K[M(X,r)]=\GK K[M(X,r)]=\rk M(X,r)=1,\] where $\rk M(X,r)$ denotes the rank of the semigroup
$M(X,r)$, i.e., the rank of the largest free commutative subsemigroup of $M(X,r)$.

On the other hand, the structure algebras of left non-degenerate finite idempotent solutions also are
right Noetherian if and only if the structure monoid is cancellative.

\begin{thm} \label{algebrastructure}
    Assume $(X,r)$ is a finite left non-degenerate idempotent solution of the Yang--Baxter equation.
    If $M=M(X,r)$ and $K$ is a field, then the following conditions are equivalent:
    \begin{enumerate}
        \item $|\Lambda|=1$.
        \item the monoid $M$ is cancellative.
        \item there exists an operation $\cdot$ on $X$ so that $(X,\cdotp)$ is a group and there exists
        a homomorphism $\varphi$ on this group so that $\lambda_x(y)=x\cdot\varphi(y)$.
        \item the algebra $K[M]$ is right Noetherian.
        \item the algebra $K[M]$ is not central (that is, its center is strictly bigger than $K$).
    \end{enumerate}
    Furthermore, if $\ch K$ does not divide $|X|$ these conditions are equivalent with
    \begin{enumerate}
        \item[(6)] the algebra $K[M]$ is semiprime.
    \end{enumerate}
    For arbitrary fields $K$, $(6)$ implies $(1)$. Moreover, if any of these conditions is satisfied
    then $K[M]$ is a subalgebra of $K[G_u]\cong K[T(G_u)][t^{\pm 1};\sigma]$, the skew Laurent polynomial
    algebra in the variable $t$ over the group algebra $K[T(G_u)]$ of the finite group $T(G_u)$.
    \begin{proof}
        The equivalence of (1) and (2) follows from Lemma~\ref{quotientgroup}.
        
        Hence, from Theorem~\ref{solgeneral}, (2) is equivalent with $X$ having a group structure for some binary
        operation~$\cdot$ so that there exists a permutation $\varphi\in\Sym(X)$ that satisfies the equations
        \eqref{fineq1}--\eqref{fineq4} (note that all $\varphi_x=\varphi$). It is readily verified that these
        equations reduce to $\varphi$ being an automorphism of the group $(X,\cdotp)$. This shows (3) holds.
        Conversely, if (3) holds then $\lambda_{x}^{-1}(x)= \varphi^{-1}(x^{-1} \cdot x)=1$, the identity of
        the group $(X,\cdotp)$. Hence, $|\Lambda |=1$ and thus (1) holds.
        
        Assume (1) holds, say $\Lambda =\{ u\}$. If $S=S(X,r)$ then Lemma~\ref{quotientgroup} yields that
        $S=S_u= \free{c_u}\circ F_u$, and $F_u$ is a finite set. Furthermore, since $\Lambda $ is a singleton,
        the element $c_u$ is central in $S_u$. Hence $K[M]=K+K[S]=\sum_{f\in F_u}f\circ K[\free{c_u}]$,
        a finite right (and left) module over a polynomial algebra in one variable. Hence, $K[M]$ is right Noetherian.
        So (4) holds.
        
        Assume (4) holds. We need to prove that $\Lambda$ is a singleton. We prove this by contradiction.
        So suppose $|\Lambda|>1$, i.e., there exist distinct $u,v\in\Lambda$ and $s_x\coloneqq(x,\lambda_x)\in S_u$
        and $s_y\coloneqq(y,\lambda_y)\in S_v$. Then, by Proposition~\ref{decomposition}, $s_x^d \circ m =s_y^d \circ m$
        for any $1\ne m\in M$. Consequently, $(s_x^d-s_y^d)\circ K[M]=(s_x^d-s_y^d)K$, and
		\[I\coloneqq K[M]\circ(s_x^d-s_y^d)\circ K[M]=K[M]\circ(s_x^d-s_y^d).\]
        We now show that the principal right ideal $I$ is not a finite right $K[M]$-module. Indeed, it were
        $I=\sum_{i=1}^nf_i\circ K[M]$ for some $f_1,\dotsc,f_n\in I$ then $I=\sum_{i=1}^nf_iK$ would be a finite-dimensional
        ideal. But the infinite set $\{(s_x^d)^j-(s_y^d)^j:j\ge 1\}$ contained in $I$ is linearly independent over $K$.
        This proves that indeed $\Lambda$ is a singleton.
        
        Also notice that if $\Lambda$ is a singleton then $c_u\in Z(K[M])$ and thus $Z(K[M])\ne K$ and (5) holds.
        To prove the converse, we need to show that $|\Lambda|>1$ implies $Z(K[M])=K$. Since $Z(K[M])$ is $\mb{N}$-graded,
        it is enough to show that there are no non-zero central homogeneous elements of $K[M]$ of positive degree.
        Suppose, on the contrary, that $0\ne\alpha=\sum_{x\in X}\alpha_x(nx,\lambda_{nx})\in Z(K[M])$ is a homogeneous
        element of degree $n>0$. Let $y,z\in X$ be such that $(y,\lambda_y)\in S_u$ and $(z,\lambda_z)\in S_v$ with $u\ne v$.
        Then $(s_y^d-s_z^d)\circ\alpha=0$ implies that $\alpha\circ(s_y^d-s_z^d)=0$. Hence $\beta=\alpha\circ s_y^d=\alpha\circ s_z^d$.
        Moreover, \[\beta=\alpha\circ s_y^d=s_y^d\circ\alpha=\sum_{x\in X}\alpha_x(du,\id)\circ(nx,\lambda_{nx})
        =\sum_{x\in X}\alpha_x((n+d)x,\lambda_{(n+d)x}).\]
        So also $\beta\ne 0$. Define $\mc{H}=\{\lambda_{nx}:x\in X\}\s\mc{G}(X,r)$ and $X_\sigma=\{x\in X:\lambda_{nx}=\sigma\}$
        for $\sigma\in \mc{H}$. Then
        \begin{align*}
            \beta=\alpha\circ s_y^d & =\sum_{x\in X}\alpha_x(nx,\lambda_{nx})\circ(du,\id)\\
            & =\sum_{x\in X}\alpha_x (\lambda_{nx} ((n+d)u),\lambda_{nx})\\
            & =\sum_{\sigma\in \mc{H}}\Big(\sum_{x\in X_\sigma}\alpha_x\Big)(\sigma((n+d)u),\sigma).
        \end{align*}
        Thus $\beta\ne 0$ yields that $(\sigma( (n+d)u,\sigma)$ lies in the support of $\beta$ for some $\sigma\in\mc{H}$,
        i.e., $0\ne \sum_{x\in X_\sigma}\alpha_x$. Because also
        \[\beta=\alpha\circ s_z^d=\sum_{\tau\in\mc{H}}\Big(\sum_{x\in X_\tau}\alpha_x\Big)(\tau((n+d)v),\tau),\]
        we get $(\sigma((n+d)u),\sigma)=(\tau((n+d)v),\tau)$ for some $\tau\in \mc{H}$. Therefore, we have $\sigma=\tau$
        and $\sigma((n+d)u)=\tau ((n+d)v)$, which leads to $\sigma(u)=\sigma(v)$, and consequently to $u=v$, a contradiction.
        This finishes the proof of (5) implies (1).
        
        To prove (6), assume $\ch K$ does not divide $|X|$. If (1) holds then $S=S_u$ is cancellative and hence so is $M$.
        By Lemma~\ref{quotientgroup}, $S$, and thus $M$ has a group of quotients $G_u\cong T(G_u)\rtimes_{\sigma_u} \mb{Z}$,
        that is a central localisation with respect the central monoid $\free{c_u}$, which is also a regular subset. It is
        well-known that $K[G_u]$ is a semiprime algebra (see, e.g., \cite[Theorem 4.2.13]{MR798076}) because $\ch K$ does
        not divide $|T(G_u)|=|X|$. Hence also $K[M]$ is semiprime. For the converse, it follows from the proof of the
        equivalence of (4) and (1) that if $|\Lambda|\ne 1$ that the ideal generated by the elements $s_x^{d}-s_y^{d}$
        is a nilpotent ideal of nilpotency class $2$. So the algebra $K[M]$ is not semiprime.
        
        The last part of the result is an immediate consequence of Lemma~\ref{quotientgroup}.
    \end{proof}
\end{thm}

Now, we give two remarks about the algebraic properties of the structure algebra determined
by a finite left non-degenerate idempotent solution.

\begin{rem}
    The structure algebra $K[M(X,r)]$ of a finite left non-degenerate idempotent solution $(X,r)$
    is prime only if $|X|=1$. Indeed if $K[M(X,r)]$ is prime then Theorem~\ref{algebrastructure}
    implies $|\Lambda|=1$. Hence, by Theorem~\ref{algebrastructure} and because $K[G_u]$ is a central
    localisation of $K[S_u]$, it follows that $K[M(X,r)]$ is prime if and only if $K[G_u]$ is prime.
    The latter only holds if $|G_u|=1$. So indeed $|X|=1$. 
\end{rem}

\begin{rem}
    Since each $K$-algebra which is a finite module over a commutative finitely generated $K$-subalgebra
    is automaton (see \cite{COAut}), it follows by Proposition~\ref{stralg} that $K[M(X,r)]$ is an automaton
    algebra in case $(X,r)$ is a finite left non-degenerate idempotent solution. It also is known
    (see, e.g., \cite{Ufna}) that each finitely generated $K$-algebra admitting a finite Gr\"obner basis
    is automaton. Because there are numerous results indicating that computational properties of such algebras
    as well as several of their algebraic and structural properties behave better than arbitrary finitely
    generated (or even finitely presented) algebras (see, e.g., \cite{MR2497577,Okn1991,MR2450723,Ufna,BBL}),
    it would be of interest to decide whether $K[M(X,r)]$ always admit a finite Gr\"obner basis.
\end{rem}

In case all $\lambda_x$ are equal we show in the following example that $K[M(X,r)]$ indeed has a finite Gr\"obner basis.
\begin{ex}
    Let $(X,r)$ be a finite idempotent solution defined as $r(x,y)=(\lambda(y),y)$, where $\lambda\colon X\to X$
    is an arbitrary map. (Note that if $(X,r)$ is an arbitrary finite left non-degenerate idempotent solution,
    written as $r(x,y)=(\lambda_x(y),\rho_y(x))$, such that $\lambda_x=\lambda_y$ for all $x,y\in X$ (call this map $\lambda$)
    then $(X,r)$ is of the above form.) Let us choose (arbitrarily) a well ordering $\le$ on $X$ (e.g., $x_1<\dotsb<x_n$ provided
    $X=\{x_1,\dotsc,x_n\}$) and extend it to the degree-lexicographic ordering on the free monoid $X^*$ freely generated by $X$.
    As \[M(X,r)=\free{X\mid xy=\lambda (y)y,\; x,y\in X}=\free{X\mid xy=yy,\; x,y\in X},\] 
    we get that $K[M(X,r)]\cong K[X^*]/I$, where the ideal $I$ of $K[X^*]$ is generated by the set $G$
    consisting of all binomials of the form:
    \begin{enumerate}
        \item $yx-xx$ for $x,y\in X$ with $x<y$,
        \item $yy-xy$ for $x,y\in X$ with $x<y$,
        \item $yz-xz$ for $x,y,z\in X$ with $x<y$.
    \end{enumerate}
    Moreover, one may check that $G$ is a Gr\"obner basis of the ideal $I$ (this can be done easily by using
    Bergman's Diamond Lemma \cite{MR506890}). In particular, the algebra $K[M(X,r)]$ admits a finite Gr\"obner basis.
\end{ex}

As an application, we obtain one of the main results, Theorem 3.10, proven in \cite{StVo21}.
Recall that a left non-degenerate solution $(X,r)$ is said to be latin if for each $x\in X$ the map
$X\to X\colon y\mapsto \lambda_y (x)$ is bijective. In the terminology of the authors this is formulated
as $(X,*)$ is a quasigroup, that is the operation $*$ on $X$ defined by $x*y=\lambda_x (y)$ defines a left
and right quasigroup, i.e., for each $x\in X$, left and right multiplications by $x$ are bijective. 

\begin{cor}\label{idempotentlatin}
    Assume $(X,r)$ is a finite left non-degenerate idempotent solution of the Yang--Baxter equation.
    Then the following conditions are equivalent:
    \begin{enumerate}
        \item $(X,r)$ is a latin solution.
        \item $X$ can be equipped with a binary operation $\cdot$ so that $(X,\cdotp)$ is a group
        and there exists $\varphi\in \Aut(X,\cdotp)$ so that $\lambda_x (y) =x\cdot\varphi (y)$.
    \end{enumerate}
    \begin{proof}
        Assume $(X,r)$ is a latin solution. Let $S=S(X,r)$ and fix $y\in X$. Because of Proposition~\ref{decomposition},
        let $u\in\Lambda$ be so that $(y,\lambda_y)\in S_u$. Also in $S$ we have that
        \[(x,\lambda_x)\circ(y,\lambda_y)=(\lambda_x(y),\lambda_{\lambda_x(y)})\circ(q(\lambda_x(y)),\lambda_{q(\lambda_x(y))})\]
        and thus $(q(\lambda_x(y)),\lambda_{q(\lambda_x(y))})\in S_u$ for all $x\in X$. Because of the assumption this yields
        $(q(z),\lambda_{q(z)})\in S_u$ for all $z\in X$. Hence $q^2$ is a constant map and thus, by Lemma \ref{qq}, it follows
        that the map $q$ is constant. As $\Img(q)=\Lambda$, we get that $\Lambda$ is a singleton. By Theorem~\ref{algebrastructure}
        it follows that (2) holds.
        
        The converse is obvious.
    \end{proof}
\end{cor}

In Theorem~\ref{algebrastructure} we dealt with idempotent left non-degenerate solutions $(X,r)$ with
$q\colon X \to \Lambda$ a constant map. Next we consider the case that the mapping $q$ is bijective, i.e., $X=\Lambda$.

\begin{cor}\label{qbijective}
    Assume $(X,r)$ is a finite left non-degenerate idempotent solution of the Yang--Baxter equation.
    Then $\Lambda=X$ if and only if $r(x,y)=(\varphi(y),y)$ for some permutation $\varphi\in\Sym(X)$.
    Such a solution we will simply denote as $(X,r_{\varphi})$.
    
    In particular, the number of non-isomorphic finite left non-degenerate idempotent solutions on a set
    $X$ with bijective diagonal map equals the number of partitions of $|X|$.
    \begin{proof}
        Because of the assumption and Theorem~\ref{solgeneral}, $X$ can be equipped with a left simple band
        structure $(X,\cdotp)$ and there exist bijective maps $q\colon X\to X$ and $\varphi_x\colon X \to X$,
        for $x\in X$, satisfying \eqref{fineq1}, \eqref{fineq2}, \eqref{fineq3} and \eqref{fineq4}.
        Because $q$ is bijective and all elements of $X$ are left identities these equations become:
        \begin{align*}
            \varphi_x \varphi_y (z)&= \varphi_{\varphi_x(y)} \varphi_{q\varphi_{x}(y)}(z),\\
            \varphi_{q\varphi_x \varphi_y (z)}q\varphi_{y}(z)&=q\varphi_{\varphi_{x}(y)} \varphi_{q\varphi_x(y)}(z),\\
            \varphi_{q\varphi_{x}(y)}(z)&= \varphi_{\varphi_x \varphi_y (z)} q\varphi_y(z),\\
            q(x) &= \varphi_x^{-1}(x).
        \end{align*}
        Put $t\coloneqq \varphi_x \varphi_y (z)$. Then the second equation becomes $\varphi_{q(t)}q\varphi_x^{-1}(t)=q(t)$
        or equivalently $q\varphi_x^{-1}(t) =\varphi_{q(t)}^{-1} q(t)$. It follows that all $\varphi_x$ are equal; let us
        denote this map simply as $\varphi$. The equations then simply say that $q=\varphi^{-1}$ and thus
        $x\cdot\varphi_x(y)=\varphi(y)$ and $q(x\cdot\varphi_x(y))=q(\varphi(y))=y$. Hence indeed $r(x,y)=(\varphi(y),y)$
        for $x,y\in X$. This proves the necessity of the result. The converse is obvious.
        
        For the last part of the statement it is sufficient to notice that solutions $(X,r_{\varphi})$ and $(X,r_{\psi})$
        are isomorphic if and only if the permutations $\varphi$ and $\psi$ are conjugate.
    \end{proof}
\end{cor}

As an application we obtain another main result, Theorem 4.7, proven in \cite{StVo21}.

\begin{cor}\label{primecase}
    The finite left non-degenerate idempotent solutions of the Yang--Baxter equation of prime cardinality $p$
    are the solutions of the following type:
    \begin{enumerate}
        \item $(X,r)$ with $X$ a set of cardinality $p$, $\varphi\in\Sym(X)$ and $r(x,y)=(\varphi(y),y)$.
        \item $(X,r)$ with $X$ having a group structure $(X,\cdotp)\cong(\mb{Z}_p,+)$, $\varphi\in\Aut(X,\cdotp)$
        and $r(x,y)=(x\cdot\varphi(y),1)$, where $1$ is the identity of the group $(X,\cdotp)$.
    \end{enumerate}
    \begin{proof}
        Let $(X,r)$ be a finite left non-degenerate idempotent solution of the Yang--Baxter equation such that $|X|=p$.
        Then Lemma~\ref{infotorsioncover} yields that $|\Lambda|=1$ or $|\Lambda|=p$. The latter case is dealt with
        in Corollary~\ref{qbijective} and hence $(X,r)=(X,r_{\varphi})$ for some $\varphi\in\Sym(X)$. The former case
        is dealt with in Theorem~\ref{algebrastructure} and thus one has a group structure on $X$. Because $p$ is prime,
        we have that $(X,\cdotp) \cong (\mb{Z}_p,+)$, $\lambda_{x}(y)=x\cdot\varphi(y)$ and $r(x,y)=(x\cdot\varphi(y),1)$,
        with $\varphi\in\Aut(X,\cdotp)$.
    \end{proof}
\end{cor}

We now give an example of a finite left non-degenerate idempotent solution $(X,r)$ that is not covered by
Corollary~\ref{idempotentlatin} and Corollary~\ref{qbijective}, i.e., $\Lambda \ne X$ and $|\Lambda|>1$.
For this it is convenient to determine when all maps $\varphi_x$ in Theorem~\ref{solgeneral} are equal.
Note that this is the case if $|\Lambda|=1$.

\begin{cor}\label{allphiequal}
    Assume $(X,r)$ is a finite left non-degenerate idempotent solution of the Yang--Baxter equation
    as in Theorem~\ref{solgeneral}. If all $\varphi_x$ are equal, and we denote this map by $\varphi$,
    then the conditions \eqref{fineq1}, \eqref{fineq2}, \eqref{fineq3} and \eqref{fineq4}
    are equivalent, respectively, with:
    \begin{enumerate}
        \item $\varphi\in\Aut(X,\cdotp)$.
        \item $\varphi q=q^2$.
        \item $q=q^4$.
        \item $q(x\cdot q^{2}(x)) =q(x)$ for all $x\in X$.
    \end{enumerate}
    \begin{proof}
        That $\varphi$ is an automorphism of $(X,\cdotp)$ is equivalent with the identity \eqref{fineq1}.
        
        The identity \eqref{fineq2} is equivalent with $\varphi q(y\cdot \varphi(z))=q(x\cdot\varphi(y)\cdot\varphi\varphi(z))$
        for all $x,y,z\in X$. Equivalently, because $\varphi$ is an automorphism,
        $\varphi q(y\cdot t)=q(x\cdot\varphi(y)\cdot\varphi(t))=q(x\cdot\varphi(y\cdot t))$. Taking into account that
        the mapping $t\mapsto y\cdot t$ is bijective, $\varphi q(z)=q(x\cdot z)$ for all $x,z\in X$. If $z\in S_u$,
        that is $q(z)=u\in\Lambda$, then, because $(X,\cdotp)$ is a left cancellative simple semigroup, there exists
        $x\in X$ so that $x\cdot z =q(z)$ (note that $q(z)$ is the identity elements of $G_u$). Hence, $\varphi q (z) =qq(z)$.
        This proves part (2).
        
        The identity \eqref{fineq3} is equivalent with $q\varphi(z)=q\varphi q(y\cdot\varphi(z))$, for all $y,z\in X$.
        Because $\varphi (z) \in G_{q(\varphi (z))}$, this is equivalent with $q\varphi (z)= q \varphi q q(\varphi (z))$,
        i.e., $q =q \varphi qq$. So by part (2), this is the same as $q= q^4$.
        
        Finally, the identity \eqref{fineq4} is equivalent with $q(x\cdot \varphi q(x)) =q(x)$, and thus by (2),
        this is the same as $q(x\cdot q^{2}(x)) =q(x)$.
    \end{proof}
\end{cor}

\begin{ex}\label{examplespecial}
    Let $G$ be a finite group, with identity $1$, and $X=\mc{M}(G,1,n,J)$, a finite left cancellative
    simple semigroup, with $n$ a positive even integer. Denote the elements of $X$ as $(g,i)$, with
    $g\in G$ and $1\le i\le n$. Let $f$ be an automorphism of the group $G$. Let $A\cup B$ be a partition
    of $\{1,\dotsc,n\}$ with $|A|=|B|$ and let $t\colon B \to A$ be a bijection. Define $\theta\colon A\cup B\to A\cup B$
    by $\theta(a)=a$ for $a\in A$ and $\theta(b)=t(b)$ for $b\in B$. So $\theta^2=\theta$. Define $q\colon X \to X$ as
    $q(g,i)=(1,\theta(i))$ and $\varphi(g,i)=(f(g),\psi(i))$, where $\psi\in\Sym(A\cup B)$ is such that $\psi|_A=\id_A$.
    Then all conditions listed in Corollary~\ref{allphiequal} are satisfied and hence it yields a solution with
    $\Lambda=\Img(q)\ne X$ and $\Lambda$ is not a singleton. 
\end{ex}

\section*{Acknowledgement}

The first named author is supported by the Engineering and Physical Sciences Research Council [grant number EP/V005995/1].
The third author is supported by National Science Centre grant 2020/39/D/ST1/01852 (Poland).
The fourth author is supported by Fonds voor Wetenschappelijk Onderzoek (Flanders), via an FWO post-doctoral fellowship, grant 12ZG221N.
The fifth author is supported by Fonds voor Wetenschappelijk Onderzoek (Flanders), via an FWO Aspirant-fellowship, grant 11C9421N.

\bibliographystyle{amsplain}
\bibliography{refs}

\end{document}